\documentclass[11pt]{article}

\usepackage{amsfonts,amssymb}
\usepackage {latexsym,times}
\usepackage {amsmath}
\usepackage{amscd}
\usepackage{epsfig,color}

\newcommand{\A}{\mathcal{A}}

\newcommand{\C}{\mathbb{C}}

\newcommand{\R}{\mathbb{R}}

\newcommand{\HH}{\mathbb{H}}

\newcommand{\K}{\mathbb{K}}

\newtheorem{theorem}{Theorem}[section]
\newtheorem{corollary}[theorem]{Corollary}
\newtheorem{lemma}[theorem]{Lemma}
\newtheorem{proposition}[theorem]{Proposition}

\newtheorem{definition}[theorem]{Definition}

\newtheorem{example}[theorem]{Example}

\newtheorem{remark}[theorem]{Remark}

\input{epsf.sty}

\begin{document}

\title{Cohomology and Formal Deformations \\ of Left  Alternative Algebras }

\author{
Mohamed
Elhamdadi\\
Department of Mathematics\\
and Statistics
\\ University of South Florida\\
\small{emohamed@math.usf.edu}
\and
Abdenacer Makhlouf\\
Laboratoire de Math\'ematiques, \\Informatique et Applications
\\ Universit\'e de Haute-Alsace\\
\small{Abdenacer.Makhlouf@uha.fr}
}

\maketitle

\begin{abstract}
The purpose of this paper is to  introduce an algebraic cohomology and formal deformation theory  of left alternative algebras.  Connections to some other algebraic structures are given also.
\end{abstract}

\section{Introduction}

Deformation theory arose mainly from geometry and physics.  In the later field,  the non-commutative associative multiplication of operators in quantum mechanics is thought of as a formal associative deformation of the pointwise multiplication of the algebra of symbols of these operators.  In the sixties, Murray Gerstenhaber introduced algebraic formal deformations for associative algebras in a series of papers (see \cite{Gerst1,Gerst2,Gerst3,Gerst4}). He used formal series and showed that the theory is intimately connected to the cohomology of the algebra. The same approach was extended to several algebraic structures. Other approaches to study deformations exist, see
\cite{Fialowski86, Fialowski88, Fialowski90, Goze88,
Goze-Remm02,LaudalLNM}, see \cite{MakhloufDeform07} for a review.

In this paper we introduce a cohomology  and a formal deformation theory of left alternative algebras. We also review the connections of alternative algebras to other algebraic structures. In Section 2, we review the basic definitions and properties related to alternative algebras.  In Section 3, we discuss in particular all the links between alternative algebras and some other algebraic structures such as Moufang loops, Malcev algebras,  Jordan algebras and Yamaguti-Lie algebras called also generalized Lie triple systems. In Section 4, we introduce a cohomology theory of left alternative algebras. We compute the second cohomology group of the $2$ by $2$ matrix algebra. It is known that, as an associative algebra its second cohomology group is trivial, but we show that this is not the case as left alternative algebra.  Finally, in Section 5, we consider the formal deformation theory of left alternative algebras.

\section{Preliminaries}
Throughout this paper $\mathbb{K}$ is an algebraically closed field of
characteristic 0.
\subsection{Definitions }

\begin{definition} {\rm \cite{McCrimmon}
 A left alternative $\K$-algebra (resp. right alternative $\K$-algebra) $(\A,\mu )$  is a vector space $\A$ over $\K$ and a bilinear multiplication $\mu$ satisfying the left alternative identity, that is, for any $x, y \in \A$,
\begin{equation}\label{LeftAlternative}
\mu (x,\mu(x, y))=\mu(\mu(x, x) , y).
\end{equation}
respectively, right alternative identity, that is
\begin{equation}\label{RightAlternative}
\mu(\mu(x,y) , y)=\mu (x,\mu(y, y)).
\end{equation}
An alternative algebra is one which is both left and right alternative algebra.}
\end{definition}
\begin{lemma}{\rm
Let $\mathfrak{as}$  denotes the associator, which is a trilinear
map defined by $\mathfrak{as}(x,y,z)=\mu(\mu(x, y), z)-\mu (x,\mu
(y, z)).$  An algebra is alternative if and only if the associator
$\mathfrak{as}(x,y,z)$ is an alternating function of its arguments,
that is
$$\mathfrak{as}(x,y,z)=-\mathfrak{as}(y,x,z)=-\mathfrak{as}(x,z,y)=
-\mathfrak{as}(z,y,x)$$
}
\end{lemma}
\noindent
This lemma implies then that the following identities are satisfied
\begin{eqnarray*}
 && \mathfrak{as}(x,x,y)=0 \quad \text{(left alternativity)},\\ && \mathfrak{as}(y,x,x)=0\quad \text{(right alternativity)} \\ && \mathfrak{as}(x,y,x)=0 \quad \text{(flexibility       ).}
 \end{eqnarray*}

\noindent
By linearization, we have the following characterization of left (resp. right) alternative algebras, which will be used in the sequel.

\begin{lemma}{\rm
A pair $(\A, \mu)$ is a left alternative $\K$-algebra (resp. right alternative $\K$-algebra) if and only if the identity
\begin{equation}\label{LeftAlternativeLineariz}
\mu (x,\mu(y, z))-\mu(\mu(x,y), z)+\mu (y, \mu(x, z))-\mu(\mu(y, x), z)=0.
\end{equation}
respectively,
\begin{equation}\label{RightAlternativeLineariz}
\mu (x,\mu (y, z))-\mu (\mu (x, y), z)+\mu (x, \mu (z, y))-\mu (\mu (x, z), y)=0.
\end{equation}
holds.}
\end{lemma}

\begin{remark}
{\rm When considering multiplication as a linear map $\mu: \A \otimes \A \rightarrow\A$, the condition (\ref{LeftAlternativeLineariz}) (resp. (\ref{RightAlternativeLineariz})) may be written
  \begin{equation}\label{LeftAlternativeLineariz2}
\mu \circ (\mu \otimes {\rm id} - {\rm id} \otimes \mu)\circ ({\rm id}^{\otimes 3}+\sigma_{1})=0.
\end{equation}
respectively
\begin{equation}\label{RightAlternativeLineariz2}
\mu \circ (\mu \otimes {\rm id} - {\rm id} \otimes \mu)\circ ({\rm id}^{\otimes 3}+\sigma_{2})=0.
\end{equation}
where ${\rm id}$ stands for the identity map and $\sigma_{1}$ and $\sigma_{2}$ stands for transpositions generating the permutation group $\mathcal{S}_3$ which are extended to   trilinear maps defined by, $\sigma_{1}(x_1\otimes x_2\otimes x_3)=x_2\otimes x_1\otimes x_3$ and $\sigma_{2}(x_1\otimes x_2\otimes x_3)=x_1\otimes x_3\otimes x_2$ for all $x_1,x_2, x_3\in \A$.
}
In terms of associators, the identities (\ref{LeftAlternativeLineariz}) (resp. (\ref{RightAlternativeLineariz})) are equivalent to \begin{equation}\mathfrak{as}+\mathfrak{as}\circ\sigma_{1}=0 \quad (\text{resp. }\ \mathfrak{as}+\mathfrak{as}\circ\sigma_{2}=0.)
\end{equation}
  \end{remark}
  \begin{remark}{\rm
  The notions of subalgebra, ideal, quotient algebra are defined in the usual way.
For general theory about alternative algebras (see \cite{McCrimmon}). The alternative algebras generalize  associative algebras. Recently, in \cite{Dzumad09}, it was shown that their operad is not Koszul. The dual operad of right alternative (resp. left alternative) algebras is defined by associativity and the identity $\mu (\mu (x,y), z)+\mu (\mu (x,z), y)=0,$ (resp. $\mu (\mu (x,y), z)+\mu (\mu (y,x), z)=0$). The dual operad of alternative algebras is defined by the associativity and the identity
$$\mu (\mu (x,y), z)+\mu (\mu (y,x), z)+\mu (\mu (z,x), y)+\mu (\mu (x,z), y)+\mu (\mu (y,z), x)+\mu (\mu (z,y), x)=0.$$}
\end{remark}

\subsection{Structure theorems and Examples}
We have these following fundamental properties:
\begin{itemize}
\item \textbf{Artin's theorem.}
In an alternative algebra the subalgebra generated by any two
elements is associative. Conversely, any algebra for which this is
true is clearly alternative. It follows that expressions involving
only two variables can be written without parenthesis unambiguously
in an alternative algebra.
\item \textbf{Generalization of Artin's theorem.}
Whenever three elements $x,y,z$ in an alternative algebra associate
(i.e. $\mathfrak{as}(x,y,z) = 0$), the subalgebra generated by those
elements is associative.

\item \textbf{Corollary of Artin's theorem.}
Alternative algebras are
power-associative, that is, the subalgebra generated by a single
element is associative. The converse need not hold: the sedenions
are power-associative but not alternative.
\end{itemize}

\begin{example}
[$4$-dimensional Alternative algebras.]  {\rm According to \cite{EGG}, p 144, there are exactly two alternative but not associative algebras of dimension 4 over any field.  With respect to a basis $e_0, e_1, e_2, e_3$, one algebra is given by the following multiplication (the unspecified products are zeros) $$ e_0^2=e_0, \; e_0e_1=e_1, \;e_2e_0=e_2, \;e_2e_3=e_1, \;e_3e_0=e_3, \;e_3e_2=- e_1. $$
The other algebra is given by  $$ e_0^2=e_0, \; e_0e_2=e_2, \;e_0e_3=e_3, \;e_1e_0=e_1, \;e_2e_3=e_1, \;e_3e_2=- e_1. $$
}
\end{example}

\begin{example}[Octonions]
{\rm The octonions were discovered in 1843 by John T. Graves who called them Octaves and independently by Arthur Cayley in 1845.
The octonions algebra which is also called Cayley Octaves or Cayley algebra is an  $8$-dimensional algebra defined with respect to a basis $u,e_1,e_2,e_3,e_4,e_5,e_6,e_7$, where $u$ is the identity for the multiplication, by the following multiplication table.  The table describes multiplying the $i$th row elements  by the $j$th column elements.

\[
\begin{array}{|c|c|c|c|c|c|c|c|c|}
  \hline
   \ & u& e_1 & e_2 & e_3 & e_4 & e_5 & e_6 & e_7 \\ \hline
   u& u& e_1 & e_2 & e_3 & e_4 & e_5 & e_6 & e_7 \\ \hline
   e_1 &e_1 & -u & e_4 & e_7 & -e_2 & e_6 & -e_5 &- e_3 \\ \hline
   e_2 &e_2 & -e_4 & -u & e_5 & e_1 & -e_3 & e_7 & -e_6 \\ \hline
   e_3 &e_3 & -e_7 & -e_5 & -u & e_6 & e_2 & -e_4 & e_1 \\ \hline
   e_4 &e_4 & e_2 & -e_1 & -e_6 & -u & e_7 & e_3 & -e_5 \\ \hline
   e_5 &e_5 & -e_6 & e_3 & -e_2 & -e_7 & -u & e_1 & e_4 \\ \hline
   e_6 &e_6 & e_5 & -e_7 & e_4 & -e_3 & -e_1 & -u & e_2 \\ \hline
   e_7 &e_7 & e_3 & e_6 & -e_1 & e_5 & -e_4 & -e_2 & -u \\
    \hline
\end{array}
\]
\\
\noindent
The octonion algebra is a typical example of alternative algebras.  As stated early the subalgebra generated by any two elements is associative.  In fact, the subalgebra generated by any two elements of the octonions is isomorphic to the algebra of reals $\R$, the algebra of complex numbers $\C$ or the algebra of quaternions $\HH$, all of which are associative. See \cite{Baez} for the role of the octonions in algebra, geometry and topology and see also \cite{Albuquerque} where octions are viewed as quasialgebra.
}
\end{example}

\section{Connections to other algebraic structures}

We begin by recalling some basics of {\it Moufang} loops, Moufang algebras and Malcev algebras.

\begin{definition}{\rm  \cite{Paal08}  Let $(M,*)$ be a set  with a binary operation.  It  is called a \emph{Moufang loop} if it is a quasigroup with an identity $e$ such that the binary operation satisfies one of the following equivalent identities:}
\begin{eqnarray}
 x*(y*(x*z))=((x*y)*x)*z, \label{mouf1}\\
 z*(x*(y*x))=((z*x)*y)*x, \label{mouf2}\\
 (x*y)*(z*x)=(x*(y*z))*x. \label{mouf3}
\end{eqnarray}
 \end{definition}
The typical examples include groups and the set of nonzero octonions which gives nonassociative Moufang loop.\\
As in the case of Lie group, there exists a notion of analytic Moufang loop \cite{Paal04,Paal08,Malcev}.  An analytic Moufang loop  $M$ is a real analytic manifold with the multiplication and the inverse, $g\mapsto g^{-1}$,  being  analytic mappings. The tangent space $T_e M$ is equipped with a skew-symmetric bracket $[~,~] :T_e M \times T_e M\rightarrow T_e M$ satisfying the Malcev's identity that is
\begin{equation}\label{MalcevIdentity}
[J(x,y,z),x]=J(x,y,[x,z])
\end{equation}
for any $x,y,z\in T_e M$ and  where $J$ corresponds to Jacobi's identity i.e.
$$J(x,y,z)=[x,[y,z]]+[y,[z,x]]+[z,[x,y]].
$$
\begin{definition}{\rm \cite{kerdman79}
A \emph{Malcev} $\K$-algebra  is a vector space over $\K$ and a skew-symmetric  bracket satisfying the identity (\ref{MalcevIdentity}).
}
\end{definition}
The Malcev algebras are also called Moufang-Lie algebras. We have the following fundamental Kerdman's theorem \cite{kerdman79}:

\begin{theorem}[Kerdman]
  For any real Malcev algebra there exists an analytic
Moufang loop whose tangent algebra is the given Malcev algebra.
\end{theorem}
 The connection to alternative algebras is given by the following proposition:

\begin{proposition}
The alternative algebras are Malcev-admissible algebras, that is the commutators define a Malcev algebra.
\end{proposition}

\noindent

\begin{remark}
{\rm
Let $\A$ be an alternative algebra with a unit. The set $U(\A)$ of all invertible elements of $\A$ forms a Moufang loop with respect to the multiplication  \cite{Shestakov03}.
 Conversely, not any Moufang loop can be imbedded into a loop of type $U(\A)$ for a suitable unital alternative algebra $\A$.  A counter-example was given in \cite{Shestakov03}. In \cite{Sandu08}, the author characterizes the Moufang loops which are imbeddable into a loop of type $U(\A)$.
 }
\end{remark}

\noindent
The Moufang algebras which
 are the corresponding algebras of a Moufang loop are defined as follows:

 \begin{definition}{\rm
 A left Moufang algebra $(\A, \mu)$ is one which is left alternative and satisfying the Moufang identity that is
 \begin{equation}
\mu( \mu(x, y), \mu(z, x))=\mu(\mu(x,\mu(y, z)), x).
 \end{equation}
 }
\end{definition}
\noindent
The Moufang identities (\ref{mouf1}, \ref{mouf2}, \ref{mouf3}) are expressed in terms of associator by
\begin{eqnarray}
\mathfrak{as}(x,y,z\cdot x)&=&x\cdot \mathfrak{as}(y,z,x)\\
\mathfrak{as}(x\cdot y,z, x)&=& \mathfrak{as}(x,y,z)\cdot x\\
\mathfrak{as}(y,x^2,z)&=&x\cdot
\mathfrak{as}(y,x,z)+\mathfrak{as}(y,x,z)\cdot x
\end{eqnarray}
\noindent
It turns out that in characteristic different from 2, all left alternative algebras are left Moufang algebras. Also, a left Moufang algebra is alternative if and only if it is flexible, that is  $\mathfrak{as}(x,y,x)=0$ for all $x,y\in \A$.

\noindent
In \cite{yamaguti1962} it shown  that  Malcev algebras form a class of a so-called  General Lie triple system, called also Lie-Yamaguti algebras (and contains Lie triple system).

\begin{definition} {\rm \cite{yamaguti1962}
 A Lie-Yamaguti algebra or General Lie triple system is an algebra $\mathcal{A}$ over a field $\K$ with a $\K$-trilinear map denoted $\tau(-,-,-)$ satisfying the following conditions:

\begin{eqnarray}
& \mu(x,x)=0, \label{}\\
&\tau (x,x,y)=0, \label{}\\
& \tau(x,y,z)+\tau(y,z,x)+\tau(z,x,y)+\mu(\mu(x,y),z)+ \mu(\mu(y,z),x)+\mu(\mu(z,x),y)=0 . \label{}\\
 &\tau(\mu(x,y),z,w)+\tau(\mu(y,z),x,w)+\tau(\mu(z,x),y,w)=0\\
& \tau(x,y,\mu(z,w))=\mu(\tau(x,y,z),w)+\mu(z,\tau(x,y,w))\\
 &\tau(x,y,\tau(z,v,w))=\tau(\tau(x,y,z),v,w)+\tau(z,\tau(x,y,v),w)+\tau(z,v,\tau(x,y,w)).
\end{eqnarray}
}
\end{definition}

\noindent
Any Lie algebra with Jacobi bracket $[-,-]$ can be Lie-Yamaguti algebra by puting $\mu(x,y):=[x,y]$ and $\tau(x,y,z):=[[x,y],z]$.  Clearly, if  $\tau(x,y,z)=0$ for all $x,y,z \in \mathcal{A}$, then the  Lie-Yamaguti algebra reduces to a Lie algebra.  And if $\mu(x,y)=0$ for all $x,y\in \mathcal{A}$, then the  Lie-Yamaguti algebra reduces to a Lie triple system (see \cite{J} for the definition of Lie triple system).

\noindent
The alternative algebras are connected to Jordan algebras as follows. Given an alternative algebra $(\A, \mu)$
then $(\A, \mu^+)$, where $\mu^+(x,y)=\mu(x,y)+\mu(y,x)$, is a Jordan
algebra,
that is the commutative multiplication $\mu^+$ satisfies the identity  $\mathfrak{as}_{\mu^+}(x^2,y,x)=0$.
\section{Cohomology of left alternative algebras}
Let $\A$ be a left alternative $\K$-algebra
defined by a multiplication $\mu $. A left alternative $p$-cochain is a linear map from $%
\A^{\otimes p}$ to $\A$. We denote by $\mathcal{C}^p ( \mathcal{A},
\mathcal{A} )$ the group of
all $p$-cochains.

\subsection{First differential map}

Let ${\rm id}$ denotes the identity map on  $\mathcal{A}$.  For  $f \in \mathcal{C}^1 ( \mathcal{A},
\mathcal{A} )$, we define the first differential $\delta ^1f \in  \mathcal{C}^2 ( \mathcal{A},
\mathcal{A} )$ by
\begin{equation}
\delta ^1f =\mu \circ \left(
f\otimes  {\rm id} \right) +\mu \circ \left( {\rm id}  \otimes f \right)
-f\circ \mu.
\end{equation}
We remark that the first differential of a left alternative algebra is similar to the first differential  map of Hochschild cohomology of an associative algebra ($1$-cocycles are derivations).

\subsection{Second differential map}
Let $\varphi\in\mathcal{C}^2 ( \mathcal{A},
\mathcal{A} )$, we define the second differential $\delta^2 \phi  \in \mathcal{C}^3 ( \mathcal{A},
\mathcal{A} )$ by,
 \begin{equation}\label{delta2}
\delta^2\phi=[\mu\circ(\phi \otimes {\rm id} - {\rm id} \otimes \phi)+\phi\circ(\mu \otimes {\rm id}-{\rm id} \otimes \mu)] \circ ({\rm id}^{\otimes 3}+\sigma_1).
\end{equation}
where $\sigma_1$ is defined on $\mathcal{A}^{\otimes 3}$ by   $\sigma_1(x\otimes y \otimes z)=y \otimes x \otimes z$.
\begin{remark}
{\rm
The left alternative algebra 2-differential defined in (\ref{delta2}) may be written using the Hochschild differential $\delta_H^2$
as
\begin{equation}
\delta^2\phi=\delta_H^2\phi\circ ({\rm id}^{\otimes 3}+\sigma_1)
\end{equation}
}
\end{remark}
\begin{proposition}
The composite $\delta ^{2}\circ \delta ^{1}$ is zero.
\end{proposition}
{\it Proof.}
Let $x, y, z \in \A$ and  $f\in \mathcal{C}^{1}({\A ,\A })$,
\begin{equation*}
\delta ^{1}f(x\otimes y) =\mu(f(x)\otimes y)+ \mu(x \otimes f(y)) -f(\mu(x\otimes y)).
\end{equation*}
Then
\begin{align*}
\delta^{2}(\delta ^{1}f )(x \otimes y \otimes z)
&= (xy)f(z)-f((xy)z)+[f(xy)]z +[xf(y)]z - [f(xy)]z +[(f(x))y]z +\\
&+(yx)f(z)-f((yx)z)+[f(yx)]z +[yf(x)]z -[f(yx)]z +[(f(y))x]z +\\
&-\{xf(yz) -f(x[yz])+[f(x)](yz)+x(yf(z))-xf(yz)+x[(f(y))z]+\\
&+yf(xz)-f(y[xz])+[f(y)](xz)+y(xf(z))-yf(xz)+y([f(x)]z) \}\\
&=[(xy)f(z)+(yx)f(z)-x(yf(z))-y(xf(z))]  - [f((xy)z)+f((yx)z)+\\
&-f(x(yz))-f(y(xz))] + [(xf(y))z+(f(y)x)z-(f(y))(xz)-x(f(y)z)]+\\
&+ [(f(x)y)z+(yf(x))z-(f(x))(yz)-y(f(x)z)] \\
&=0.
\end{align*}
After simplifying the terms which cancel in pairs, we group the remaining ones into brackets so each bracket cancels using the left alternative algebra axiom (equation (\ref{LeftAlternativeLineariz})).
$\Box$

\begin{example}\label{lemmaM22}{\rm
Let $\mathcal{A}=\mathcal{M}_{2}(\K)$ denotes the associative algebra of $2$ by $2$ matrices over the field $\K$, considered as left alternative algebra of dimension 4.  Let $e_1, e_2, e_3$ and $e_4$ be a basis of $\mathcal{A}$.  The second cohomology $H^2(\mathcal{A},\mathcal{A})$ is three-dimensional generated by $[f_1], [f_2]$ and $[f_3]$ where
\begin{eqnarray*}
f_1(e_2\otimes e_4)&=&e_1, \; f_1(e_3\otimes e_2)=-e_3, \; f_1(e_4\otimes e_1)=e_3,\;  f_1(e_4\otimes e_2)=e_4, \\
f_2(e_2\otimes e_3)&=&e_2, \; f_2(e_3\otimes e_1)=-e_4, \; f_2(e_3\otimes e_3)=e_3,\;  f_2(e_3\otimes e_4)=e_4, \\
f_3(e_2\otimes e_3)&=&e_1,\; f_3(e_3\otimes e_2)=e_4.\\
\end{eqnarray*}
The non-specified terms of these generators are zeros.  These generators were obtained independently using the softwares Maple and Mathematica.}
\end{example}

\subsection{Third differential map}
Let $\psi\in\mathcal{C}^3 ( \mathcal{A},
\mathcal{A} )$,  we define the third differential $\delta^3\psi \in \mathcal{C}^4 ( \mathcal{A},
\mathcal{A} )$ as,
\begin{eqnarray*}
\delta^3\psi  &=&  \mu( \psi\otimes {\rm id} )
( {\rm id}^{\otimes 3}-\sigma_1)+\mu( {\rm id} \otimes \psi)
( id^{\otimes 3}-\sigma_2)\\
& & -\psi(\mu\otimes {\rm id}^{\otimes 2} )
( {\rm id}^{\otimes 3}+\sigma_2\circ\sigma_1)
 +\psi ({\rm id} \otimes \mu\otimes {\rm id} )
 ( {\rm id}^{\otimes 3}+\sigma_1\circ\sigma_2)
 -\psi( {\rm id}^{\otimes 2}\otimes \mu )( {\rm id}^{\otimes
3}-\sigma_1).
\end{eqnarray*}

That is for all  $\psi\in\mathcal{C}^3 ( \mathcal{A},
\mathcal{A} )$ and  $x_1,\ldots,x_4\in \A$
\begin{eqnarray*}
\delta^3\psi (x_1,x_2,x_3,x_4) &=&
\mu(x_1 \otimes \psi(x_2\otimes x_3\otimes x_4))-
\mu(x_1 \otimes \psi(x_3\otimes x_2\otimes x_4)) \\
& &
+\mu( \psi(x_1\otimes x_2\otimes x_3)\otimes x_4 )-
\mu( \psi(x_2\otimes x_1\otimes x_3) \otimes x_4 ) \\
& &
-\psi(\mu( x_1\otimes x_2)\otimes x_3\otimes x_4 )-
\psi(\mu( x_2\otimes  x_3)\otimes x_1 \otimes x_4 )\\
& &
 +\psi (x_1\otimes \mu(x_2\otimes x_3)\otimes x_4 )+
 \psi(x_3\otimes \mu(x_1\otimes x_2) \otimes x_4 )\\
& &
-\psi( x_1\otimes x_2\otimes \mu(x_3\otimes x_4) )+
\psi( x_2\otimes x_1\otimes \mu(x_3 \otimes x_4) ).
\end{eqnarray*}

\begin{proposition}
The composite $\delta ^{3}\circ \delta ^{2}$ is zero.
\end{proposition}
\noindent
{\it Proof.}
Let $x_1,\ldots,x_4\in \A$ and  $f\in \mathcal{C}^{2}({\A ,\A })$,
\
Then, by substituting $\psi$ with $\delta ^{2}f$ in the previous formula and rearranging the terms we get
\begin{eqnarray*}
\delta^{3}(\delta ^{2}f )(x_1\otimes x_2\otimes x_3\otimes x_4)&=&
x_1 [ \delta ^{2}f(x_2\otimes x_3\otimes x_4)- \delta ^{2}f(x_3\otimes x_2\otimes x_4)] \\
& &
-[\delta ^{2}f(x_1 x_2\otimes x_3\otimes x_4 )  -\delta ^{2}f(x_3\otimes x_1 x_2 \otimes x_4 )  ]   \\
& &
 +[\delta ^{2}f (x_1\otimes x_2 x_3\otimes x_4 )-
\delta ^{2}f(x_2 x_3\otimes x_1 \otimes x_4 )]
 \\
& &
-[\delta ^{2}f( x_1\otimes x_2\otimes x_3 x_4 )-
\delta ^{2}f( x_2\otimes x_1\otimes x_3  x_4) ]\\
& &
+[\delta ^{2}f(x_1\otimes x_2\otimes x_3)- \delta ^{2}f(x_2\otimes x_1\otimes x_3) ]x_4  \\
&= &0,
\end{eqnarray*}
since $\delta ^{2}f (x \otimes y \otimes z)=\delta ^{2}f (y \otimes x \otimes z)$,  for all $x,y,z \in \A$.
$\Box$\\

It is an interesting problem to find the higher $p^{th}$ differential maps
and study the properties of the cohomology groups. This will be considered by the authors in a forthcoming work.

\noindent

%
%
%
%
%
%
%

\section{Formal Deformations of  left alternative algebras}

Let  $({\A },\mu _{0})$ be a left alternative
algebra. Let $\K [[t]]$ be the power series ring in one
variable $t$ and coefficients in $\K $ and $\A[[t]]$ be the set of
formal power series whose coefficients are elements of $\A$
(note that $\A[[t]]$ is obtained   by extending the coefficients domain of $\A$
from $\K $ to $\K [[ t]]$). Then $\A[[t]]$ is a $\K[[t]]$-module.
When $\A$ is finite-dimensional, we have $ \A[[t]]=\A\otimes _{\K
}\K[[t]]$. One notes that $V$ is a submodule of $\A[[t]]$. Given a
 $\K$-bilinear map $f :\A\times \A  \rightarrow \A$, it admits naturally an
 extension to a $\K[[t]]$-bilinear map
 $f :\A[[t]]\otimes \A[[t]]  \rightarrow \A[[t]]$, that is,
 if  $x=\sum_{i\geq0}{a_i t^i}$ and $y=\sum_{j\geq0}{b_j t^j}$ then
$f(x\otimes y)=\sum_{i\geq0,j\geq0}{t^{i+j}f (a_i\otimes b_j)}$.

\begin{definition}{\rm
Let $({\A },\mu _{0})$ be a left alternative  algebra. A \emph{formal  left alternative deformation} of
${\A }$  is given by the $\K[[t]]$-bilinear
 map $\mu_{t} :\A[[t]]\otimes \A[[t]]  \rightarrow
\A[[t]]$  of the form
$\mu_{t} =\sum_{i\geq 0}\mu_{i}t^{i},
$
where each $\mu_{i}$ is a $\K$-bilinear map $\mu_{i}:  \A\otimes \A
\rightarrow  \A$ (extended to be $\K[[t]]$-bilinear), such that
for $x, y,z\in
 \A$, the following
formal  left alternativity condition holds
\begin{equation}
\label{equ1} \mu_{t}(x\otimes \mu_{t}(y\otimes z))-\mu_{t}(\mu_{t}(x\otimes y)\otimes z))+\mu _t\left( y\otimes \mu _t\left( x\otimes z\right) \right)-\mu _t\left( \mu _t\left( y\otimes x\right)
\otimes z\right) =0.
\end{equation}
}
\end{definition}

\subsection{Deformation equation and Obstructions}
The first problem is to give
conditions about $\mu _i$ such that the deformation $\mu _t$ be
alternative.
Expanding the left side of the equation (\ref{equ1}) and collecting the coefficients of $%
t^k$ yields%
\begin{eqnarray*}
\left\{ \sum_ { i+j=k \ \ i,j\geq 0}{  \mu _i\left( x\otimes\mu _j\left( y\otimes z\right)
\right)-\mu _i\left( \mu _j\left(
x\otimes y\right)\otimes z\right)+\mu _i\left( y\otimes \mu _j\left( x\otimes z\right)
\right)-
\mu _i\left( \mu _j\left(
y\otimes x\right) \otimes z\right) } =0,\right.\\ k=0,1,2,\ldots
\end{eqnarray*}

\noindent
This infinite system, called the {\it deformation equation,} gives the
necessary and sufficient conditions for the left alternativity of $\mu _t$.
 It may be written
\begin{eqnarray}\label{equaDefo}
\left\{ \sum_{i=0}^k{\mu _i\left( x\otimes \mu _{k-i}\left(
y\otimes z\right) \right)-\mu _i\left( \mu
_{k-i}\left( x,y\right) \otimes z\right) + \mu _i\left( y\otimes \mu _{k-i}\left(
x\otimes z\right) \right) -\mu _i\left( \mu
_{k-i}\left( y\otimes x\right) \otimes z\right)}=0,\right.\\  k=0,1,2,\ldots  \nonumber
\end{eqnarray}
\noindent
The first equation $\left( k=0\right) $ is the left alternativity condition for $%
\mu _0.$\\The second shows that $\mu _1$ must be a 2-cocycle for the
Alternative algebra cohomology defined above $\left( \mu _1\in Z^2\left( \A,\A\right) \right) $.\\%
More generally, suppose that $\mu _p$ be the first non-zero coefficient
after $\mu _0$ in the deformation $\mu _t$. This $\mu _p$ is called the {\it %
infinitesimal} of $\mu _t$.

\begin{theorem} The map $\mu _p$ {\it is a} 2-{\it cocycle of the left alternative algebras
cohomology of }$\A$ {\it with coefficient in itself.}
\end{theorem}
\noindent
{\it Proof.} In the equation (\ref{equaDefo}) make the following substitution $k=p$ and $\mu
_1=\cdots =\mu _{p-1}=0$.
$\Box$

\begin{definition}{\rm The 2-cocycle $\mu _p$ is said integrable if it
is the first non-zero term, after $\mu _0,$ of a left alternative deformation.
}
\end{definition}
The integrability of $\mu _p$ implies an infinite sequence of relations
which may be interpreted as the vanishing of the obstruction to the
integration of $\mu _p$.

\noindent
For an arbitrary $k$, with $k>1,$ the $k^{th}$ equation of the system (\ref{equaDefo}) may be written%
\begin{eqnarray*}
&&
\delta ^2\mu _k\left( x\otimes y\otimes z\right) = \\&&  \sum_{i=1}^{k-1}\mu _i\left( \mu
_{k-i}\left( x\otimes y\right) \otimes z\right) -\mu _i\left( x \otimes \mu _{k-i}\left(
y\otimes z\right) \right)+\mu _i\left( \mu
_{k-i}\left( y\otimes x\right) \otimes z\right) -\mu _i\left( y\otimes \mu _{k-i}\left(
x\otimes z\right) \right).
\end{eqnarray*}
Suppose that the truncated deformation $\mu _t=\mu _0+t\mu _1+t^2\mu
_2+\cdots +t^{m-1}\mu _{m-1}$ satisfies the deformation equation. The truncated
deformation is extended to a deformation of order $m$, i.e.  $\mu _t=\mu
_0+t\mu _1+t^2\mu _2+ \cdots +t^{m-1}\mu _{m-1}+t^m\mu _m$ satisfying the
deformation equation if
\begin{eqnarray*}
&& \delta ^2\mu _m\left( x\otimes y\otimes z\right) =\\&&  \sum_{i=1}^{m-1}\mu _i\left( \mu
_{m-i}\left( x\otimes y\right) \otimes z\right) -\mu _i\left( x\otimes \mu _{m-i}\left(
y\otimes z\right) \right)+\mu _i\left( \mu
_{m-i}\left( y\otimes x\right) \otimes z\right) -\mu _i\left( y\otimes \mu _{m-i}\left(
x\otimes z\right) \right).
\end{eqnarray*}
The right-hand side of this equation is called the {\it obstruction }to
finding $\mu _m$ extending the deformation.\\\\
\noindent
We define  a  square operation on $2$-cochains by
$$\mu _i\Box \mu _j\left( x\otimes y\otimes z\right) =\mu _i\left( \mu _j\left(
x\otimes y\right) \otimes z\right) -\mu _i\left( x\otimes \mu _j\left( y\otimes z\right) \right)+
\mu _i\left( \mu _j\left(
y\otimes x\right)\otimes z\right) -\mu _i\left( y\otimes \mu _j\left( x\otimes z\right) \right)
,$$ then  the obstruction may be written
$
\sum_{i=1}^{m-1}{\mu_i\Box \mu _{m-i}}\text{ or }
\sum_{ i+j=m \ \ i,j\neq m}{ \mu _i\Box \mu _j }.
$\\

\noindent
A straightforward computation gives the following
\begin{theorem} The obstructions are left alternative 3-cocycles.
\end{theorem}

%
\begin{remark}
{\rm
\begin{enumerate}
\item
 The cohomology class of the element $\sum_ {i+j=m,\ \ i,j\neq m}
\mu _i\Box \mu _j$ is the first obstruction to the integrability of $\mu _m$.%
\\Let us consider now how to extend an infinitesimal deformation to a deformation of order 2. Suppose $m=2$ and $\mu _t=\mu _0+t\mu _1+t^2\mu _2$. The deformation
equation of the truncated deformation of order 2 is equivalent to the finite system:%
$$
\left\{
\begin{array}{lll}
\mu _0\Box \mu _0&=&0\quad \left( \mu _0
\text{ is left alternative}\right)  \\ \delta \mu _1&=&0\quad \left( \mu _1\in
Z^2\left( \A,\A\right) \right)  \\
\mu _1\Box \mu _1&=&\delta \mu _2
\end{array}
\right.
$$

Then $\mu _1\Box \mu _1$ is the first obstruction to integrate $\mu _1$ and
$\mu _1\Box \mu _1\in Z^3\left( \A,\A\right) .$\\The elements $\mu _1\Box
\mu _1$ which are coboundaries permit to extend the deformation of order one
to a deformation of order 2. But the elements of $H^3\left( \A,\A\right) $
gives the obstruction to the integrations of $\mu _1$.

\item  If $\mu _m$ is integrable then the cohomological class of $\sum_
{i+j=m, \quad i,j\neq m} \mu _i\Box \mu_j$ must be 0.\\In the previous
example $\mu _1$ is integrable implies $\mu _1\Box \mu _1=\delta \mu _2$
which means that the cohomology class of $\mu _1\Box \mu _1$ vanishes.

\end{enumerate}
}
\end{remark}
\begin{corollary}
If $H^3\left( \A,\A\right) =0$ then all obstructions vanish and every $%
\mu _m\in Z^2\left( \A,\A\right) $ is integrable.
\end{corollary}

\subsection{Equivalent and trivial deformations}

In this section, we characterize  equivalent as well as trivial
deformations of left alternative algebras.

\begin{definition}{\rm
Let  $(\A,\mu_0)$ be a left alternative algebra   and let $(\A_t,\mu_t)$ and
$(\A'_t,\mu'_t)$ be two left alternative
 deformations of $\A$,  where
 $\mu_{t}=\sum_{i\geq 0}t^{i}\mu_{i}$ and $\mu'_{t}=\sum_{i\geq 0}t^{i}\mu'_{i}$,
 with
 $\mu_{0}=\mu'_{0}$.\\
 We say that the two deformations are \emph{equivalent} if there
exists a formal isomorphism $\Phi_{t}: \A[[t]]\rightarrow \A[[t]]$,
i.e. a $\K[[t]]$-linear map that may be written in the form $
\Phi_{t}=\sum_{i\geq 0}t^{i}\Phi _{i} ={\rm id}+t\Phi
_{1}+t^{2}\Phi_{2}+\ldots$, where $\Phi_{i}\in End_{\K }(\A)$ and
$\Phi_{0}=id $ are such that the
following relations hold
\begin{equation}
\label{equIso} \Phi_{t}\circ
\mu_{t}=\mu'_{t}\circ(\Phi _{t}\otimes \Phi _{t}
).
\end{equation}
A deformation $\A_{t}$ of $\A_{0}$ is said to be \emph{trivial} if
and only if $\A_{t}$ is equivalent to $\A_{0}$ (viewed  as a
left alternative algebra  on $\A[[t]]$).
}
\end{definition}
\noindent
We discuss in the following the equivalence of two deformations. Condition \eqref{equIso} may be written as
\begin{equation}\label{isom1}
\Phi _{t}(\mu_{t}(x\otimes y)) =\mu'_{t}(\Phi _{t}(x))\otimes \Phi
_{t}(y)),\quad\forall x,y\in \A .
\end{equation}
\noindent
Equation \eqref{isom1} is equivalent to
\begin{equation}
\sum_{i\geq 0}\Phi _{i}\left(\sum_{j\geq 0}\mu_{j} (x\otimes
y)t^{j}\right)t^{i} =\sum_{i\geq 0}\mu'_{i} \left( \sum_{j\geq
0}\Phi _{j}(x)t^{j}\otimes \sum_{k\geq 0}\Phi_{k}(y)t^{k}
 \right) t^{i}
\end{equation}
or
\begin{equation*}
\sum_{i,j\geq 0}\Phi _{i}(\mu_{j}(x\otimes y))t^{i+j} =
\sum_{i,j,k\geq 0}\mu'_{i}(\Phi _{j}(x)\otimes
\Phi_{k}(y))t^{i+j+k}.
\end{equation*}
By identification of  the coefficients, one obtains
that the constant coefficients are identical,
i.e.
\begin{equation*}
\mu_{0}=\mu'_{0} \quad\text{because}\quad \Phi
_{0} =id.
\end{equation*}
For  the coefficients of $t$ one finds
\begin{equation}
\Phi _{0}(\mu_{1}(x\otimes y))+\Phi _{1}(\mu_{0}(x\otimes y))
= \mu'_{1}(\Phi _{0}(x)\otimes \Phi _{0}(y)) +
\mu'_{0}(\Phi _{1}(x)\otimes \Phi _{0}(y)) +
\mu'_{0}(\Phi _{0}(x)\otimes \Phi _{1}(y)).
\end{equation}
Since $\Phi _{0}=id$, it follows that
\begin{equation}
\mu_{1}(x, y)+\Phi _{1}(\mu_{0}(x\otimes y)) =
\mu'_{1}(x\otimes y)+\mu_{0}(\Phi _{1}(x)\otimes y) +
\mu_{0}(x\otimes \Phi _{1}(y)).
\end{equation}
Consequently,
\begin{equation}\label{equiv1}
\mu'_{1}(x\otimes y) = \mu_{1}(x\otimes y)+\Phi_{1}(\mu_{0}(x\otimes y))-\mu_{0}(\Phi
_{1}(x)\otimes y) - \mu_{0}(x\otimes \Phi _{1}(y)).
\end{equation}

\noindent
The second order conditions of the equivalence between two
deformations of a left alternative algebra
are given by \eqref{equiv1} which may be written

\begin{equation}\label{equiv1P}
\mu'_{1}(x\otimes y) = \mu_{1}(x\otimes y)-\delta^1\Phi_{1}(x\otimes y) .
\end{equation}

\noindent
In general, if the deformations $\mu _t$ and $\mu _t^{\prime
}$ of $\mu _0$ are equivalent then $\mu _1^{\prime }=\mu _1+\delta^1
f_1$.

Therefore, we have the following proposition:

\begin{proposition}
 The integrability of $\mu _1$ depends only on its cohomology class.
 \end{proposition}

\noindent
Recall that two elements are cohomologous if their difference is a
coboundary. \\
The equation $\delta^2 \mu _1=0$ implies that $\delta^2 \mu _1^{\prime }=\delta^2
\left( \mu _1+\delta^1 f_1\right) =\delta^1 \mu _1+\delta^2 \left( \delta^1
f_1\right) =0$. \\If $\mu _1=\delta^1 g$ then $\mu _1^{\prime }=\delta^1
g-\delta^1 f_1=\delta^1 \left( g-f_1\right) .$

\noindent Then if two integrable 2-cocycles are cohomologous, then the corresponding deformations are equivalent.
 
\begin{remark}{\rm
Elements of $H^2\left( \A,\A\right) $ give the
infinitesimal deformations ( $\mu _t=\mu _0+t\mu _1).$}
\end{remark}

\begin{proposition}
Let $({\A },\mu_{0})$ be a left alternative algebra. There is, over $\K [[t]]/t^2,$
a one-to-one correspondence between the
elements of $\mathit{H^{2}}(\A  ,\A)$ and the
infinitesimal deformation of ${\A }$ defined by
\begin{equation}
\mu_{t}(x\otimes y) =\mu_{0}(x\otimes y)+t \mu_{1}(x\otimes y), \quad \forall x,y\in \A.
\end{equation}

\end{proposition}
{\it Proof.} The deformation equation is equivalent to
$\delta^2 \mu_{1}=0$, that is $\mu_{1}\in\mathit{Z^{2}}(\A
,\A)$.
$\Box$

\begin{theorem}
Let $({\A },\mu_{0})$ be a left alternative algebra and  $\mu _t$ {\it be a one parameter family of
deformation of }$\mu _0$. Then $\mu _t$ is equivalent to $\mu _t=\mu
_0+t^p\mu _p^{\prime }+t^{p+1}\mu _{p+1}^{\prime }+\cdots, $ where $\mu
_p^{\prime }\in Z^2\left( \A,\A\right) $ and $\mu _p^{\prime }\notin B^2\left(
\A,\A\right) $.
\end{theorem}
{\it Proof.}
Suppose now that $\mu _t=\mu _0+t\mu _1+t^2\mu _2+\cdots, $ is a one parameter
family of deformation of $\mu _0$ for which $\mu _1=\cdots =\mu _{m-1}=0.$
The deformation equation implies $\delta \mu _m=0$ $\left( \mu _m\in
Z^2\left( \A,\A\right) \right) $. If further $\mu _m\in B^2\left( \A,\A\right) $
(ie. $\mu _m=\delta g)$, then setting the morphism $f_t=id+tf_m$,  we have, for all $x,y\in \A, $
$$
\mu _t^{\prime }(x\otimes y)=f_t^{-1}\circ \mu _t\circ (f_t\left( x\right)
\otimes f_t\left( y\right) )=\mu _0\left( x\otimes y\right) +t^{m+1}\mu
_{m+1}\left( x\otimes y\right) \cdots.
$$
And again $\mu _{m+1}\in Z^2\left( \A,\A\right) .$
$\Box$

\begin{corollary} If $H^2\left( \A,\A\right) =0,$ then all deformations of $\A$
are equivalent to a trivial deformation.
\end{corollary}
In fact, assume that there exists a non trivial deformation of $\mu
_0$. Following the previous theorem, this deformation is equivalent
to $\mu
_t=\mu _0+t^p\mu _p^{\prime }+t^{p+1}\mu _{p+1}^{\prime }+\cdots $ where $%
\mu _p^{\prime }\in Z^2\left( \A,\A\right) $ and $\mu _p^{\prime }\notin
B^2\left( \A,\A\right) $. But this is impossible because $H^2\left( \A,\A\right) =0.$

\begin{remark}{\rm
A left alternative algebra for which every formal deformation is equivalent to a trivial deformation is called rigid. The previous corollary provide a sufficient condition for a left alternative algebra to be rigid. In general this condition is not necessary.}
\end{remark}

\noindent
{\bf Acknowledgment} M. E. would like to thank S. Carter and M. Saito for fruitfull discussions.

\end{document}